\documentclass[12pt]{amsart}
\usepackage{amsfonts}
\usepackage{amsmath}
\usepackage{amssymb}
\usepackage{mathrsfs}
\usepackage{a4}
\usepackage{amscd}

\newcommand{\heute}{25 September 2008}

\theoremstyle{plain}
\newtheorem{theorem}{Theorem}[section]

\newtheorem{lemma}[theorem]{Lemma}
\newtheorem{corollary}[theorem]{Corollary}

\theoremstyle{remark}

\newtheorem{definition}[theorem]{Definition}

\newtheorem*{defn}{Definition}
\newtheorem*{rk}{Remark}
\newtheorem*{bsp}{Example}
\newtheorem*{notn}{Notation}

\newcommand{\dashTwo}[1]{\textup{(\ref{two}${}'$)}}

\newcommand{\ignore}[1]{}

\newcommand{\f}[1][p]{\mathbb{F}_{#1}}

\newcommand{\Gro}[1]{Gr\"ob\-ner}

\newcommand{\Ess}{\operatorname{Ess}}

\newcommand{\Hom}{\operatorname{Hom}}

\newcommand{\eps}{\varepsilon}

\DeclareMathOperator{\Res}{Res}

\newcommand{\rP}{\mathcal{P}}

\DeclareMathSymbol\normal{\mathrel}{AMSa}{"43}

\newcommand{\Ste}{\mathcal{A}}

\newcommand{\GL}{\mathit{GL}}
\newcommand{\abs}[1]{\left|#1\right|}
\newcommand{\SL}{\mathit{SL}}
\newcommand{\Mui}[1]{M{\`u}i}

\begin{document}

\title[Essential cohomology]{Essential cohomology for elementary
abelian $p$-groups}
\author[F.~Altunbulak Aksu]{Fatma Altunbulak Aksu}
\address{Dept of Mathematics \\
Bilkent University \\ Bilkent, 06800 \\ Ankara, Turkey}
\email{fatma@fen.bilkent.edu.tr}
\thanks{The first author was supported by a Ph.D.\@ research
scholarship from the Scientific and Technical Research Council
of Turkey (T\"UB\.ITAK-BAYG)}
\author[D.~J. Green]{David J. Green}
\address{Dept of Mathematics \\
Friedrich-Schiller-Universit\"at Jena \\ 07737 Jena \\ Germany}
\email{David.Green@uni-jena.de}
\date{\heute}

\begin{abstract}
\noindent
For an odd prime~$p$ the cohomology ring of an elementary abelian
$p$-group is polynomial tensor exterior. We show that the
ideal of essential classes is the Steenrod closure of the
class generating the top exterior power. As a module over the
polynomial algebra, the essential ideal is free on the set of
\Mui. invariants.
\end{abstract}

\maketitle

\section{Introduction}
\label{sect:intro} \noindent Let $G$ be a finite group and $k$ a
field whose characteristic~$p$ divides the order of~$G$. A
cohomology class $x \in H^n(G,k)$ is called \emph{essential} if its
restriction $\Res_H(x)$ is zero for every proper subgroup $H$~of
$G$. The essential classes form an ideal, called the essential ideal
and denoted by $\Ess(G)$. It is standard that restriction to a Sylow
$p$-subgroup of~$G$ is a split injection (see for example Theorem
XII,10.1 of~\cite{CartanEilenberg}), and so the essential ideal can
only be nonzero if $G$ is a $p$-group. Many $p$-groups have nonzero
essential ideal, for instance the quaternion group of order eight.
The essential ideal plays an important role and has therefore been
the subject of many studies: two such being Carlson's work on the
depth of a cohomology ring~\cite{Carlson:DepthTransfer}, and the
cohomological characterization due to Adem and Karagueuzian of those
$p$-groups whose order~$p$ elements are all
central~\cite{AdKa:Ess}\@.

The nature of the essential ideal depends crucially on whether or
not the $p$-group $G$ is elementary abelian. If $G$ is not
elementary abelian, then a celebrated result of Quillen (Theorem
7.1~of \cite{Quillen:SpectrumI}) implies that $\Ess(G)$ is a
nilpotent ideal. By contrast, the essential ideal of an elementary
abelian $p$-group contains non-nilpotent classes. Work to date on
the essential ideal has concentrated on the non-elementary abelian
case. In this paper we give a complete treatment of the outstanding
elementary abelian case. As we shall recall in the next section, the
case $p=2$ is straightforward and well known. So we shall
concentrate on the case of an odd prime~$p$\@.

So let $p$ be an odd prime and $V$ a rank~$n$ elementary abelian $p$-group.
We may equally well view $V$ as an $n$-dimensional $\f$-vector space.
Recall that the cohomology ring has the form
\begin{equation}
H^*(V,\f) \cong S(V^*) \otimes_{\f} \Lambda(V^*) \, ,
\label{eqn:Vcoho}
\end{equation}
where the exterior copy of the dual space~$V^*$ is $H^1(V,\f)$, and
the polynomial copy lies in $H^2(V,\f)$: specifically, the
polynomial copy is the image of the exterior copy under the
Bockstein boundary map~$\beta$.
Our first result is as follows:

\begin{theorem}
\label{thm:Steenrod}
Let $p$~be an odd prime and $V$ a rank~$n$ elementary abelian $p$-group.
Then the essential ideal $\Ess(V)$ is the Steenrod closure
of~$\Lambda^n(V^*)$. That is, $\Ess(V)$ is the smallest ideal in
$H^*(V,\f)$ which contains the one-dimensional space
$\Lambda^n(V^*) \subseteq H^n(V,\f)$
and is closed under the action of the Steenrod algebra.
\end{theorem}

\noindent
Our second result concerns the structure of $\Ess(V)$ as a module
over the polynomial subalgebra $S(V^*)$ of $H^*(V,\f)$. It was conjectured
by Carlson (Question 5.4 in~\cite{Carlson:Problems}) -- and earlier in a less
precise form by \Mui.~\cite{Mui:Essay} -- that the essential ideal
of an arbitrary $p$-group is free and finitely generated as a module
over a certain polynomial subalgebra of the cohomology ring. In~\cite{essCM},
the second author demonstrated finite generation, and for most $p$-groups
of a given order was able to prove freeness as well: specifically
the method works provided the group is not a direct product in which
one factor is elementary abelian of rank at least two. Our second result
states that Carlson's conjecture holds for elementary abelian $p$-groups too,
and gives explicit free generators.

\begin{theorem}
\label{thm:free}
Let $p$~be an odd prime and $V$ a rank~$n$ elementary abelian $p$-group.
Then as a module over the polynomial part $S(V^*)$ of the cohomology
ring $H^*(V,\f)$, the essential ideal $\Ess(V)$ is free on the set of
\Mui. invariants, as defined in Definition~\ref{defn:Mui}\@.
\end{theorem}

\paragraph{Structure of the paper}
In \S\ref{section:mod2} we briefly cover the well-known case $p=2$.
We introduce the \Mui. invariants in~\S\ref{section:Mui}\@. After proving
Theorem~\ref{thm:free} in~\S\ref{section:jointAnn} we consider the action of
the Steenrod algebra on the \Mui. invariants in order to prove
Theorem~\ref{thm:Steenrod} in~\S\ref{section:Steenrod}\@.

\section{Elementary abelian $p$-groups and the case $p=2$}
\label{section:mod2}
\noindent
The cohomology group $H^1(G,\f)$ may be identified with the set of
group homomorphisms $\Hom(G,\f)$. This set is an $\f$-vector space, and
-- assuming that $G$ is a $p$-group --
the maximal subgroups of $G$ are in bijective correspondence
with the one-dimensional subspaces: the maximal subgroup corresponding to
$\alpha \colon G \rightarrow \f$ being $\ker(\alpha)$. Of course, the
cohomology class $\alpha \in H^1(G,\f)$ has zero restriction to the maximal
subgroup $\ker(\alpha)$. Note that in order to determine $\Ess(G)$ it
suffices to consider restrictions to maximal subgroups.

\begin{defn}
Denote by $L_n$ the polynomial
\[
L_n(X_1,\ldots,X_n) = \det \begin{vmatrix}
X_1 & X_2 & \cdots & X_n \\
X_1^p & X_2^p & \cdots & X_n^p \\
\vdots & \vdots & \ddots & \vdots \\
X_1^{p^{n-1}} & X_2^{p^{n-1}} & \cdots & X_n^{p^{n-1}} \end{vmatrix}
\in \f\lbrack X_1,\ldots,X_n \rbrack \, .
\]
\end{defn}

\noindent
There is a well-known alternative description of~$L_n$.


\begin{lemma}
\label{lemma:eqnLn}
$L_n$ is the product of all monic linear forms in $X_1,\ldots,X_n$.
So for an $n$-dimensional $\f$-vector space $V$ we
may define $L_n(V) \in S(V^*)$ up to a nonzero scalar multiple by
\begin{equation}
\label{eqn:Ln}
L_n(V) = \prod_{[x] \in \mathbb{P}V^*} x \, .
\end{equation}
\end{lemma}

\begin{proof}
First part: Here we call a linear form monic if the first nonzero coefficient
is one.
The right hand side divides the left. Both sides have the same
total degree. And the coefficient of $X_1 X_2^p X_3^{p^2} \cdots X_n^{p^{n-1}}$
is $+1$ in both cases. The second part follows.
\end{proof}

\noindent
Let $V$ be an elementary abelian $2$-group. Then $H^*(V,\f[2]) \cong S(V^*)$,
where the dual space~$V^*$ is identified with $H^1(V,\f[2])$. Pick
$x_1,\ldots,x_n$ to be a basis for $H^1(V,\f[2])$.
The following is well-known:

\begin{lemma}
\label{lemma:p2}
For an elementary abelian $2$-group~$V$, the essential ideal
is the principal ideal in $H^*(V,\f[2])$ generated by $L_n(x_1,\ldots,x_n)$.

Moreover, $\Ess(V)$ is the free $S(V^*)$-module on~$L_n(V)$,
and the Steenrod closure of this one generator.
\end{lemma}

\begin{proof}
$L_n(V)$ is essential, because every nonzero linear form is a factor
and every maximal subgroup is the kernel of a nonzero linear form.
Now suppose that $y$ is essential, and let $x \in V^*$ be a nonzero
linear form. Let $U \subseteq V^*$ be a complement of the subspace
spanned by~$x$. So $y = y'x + y''$ with $y' \in S(V^*)$ and $y'' \in
S(U)$. Hence $\Res_H (y'') = 0$ for $H = \ker(x)$, as $y$ is
essential and $\Res_H(x)=0$. But the map $\Res_H \colon V^*
\rightarrow H^*$ satisfies $\ker(\Res_H) \cap U = 0$, and so
$\Res_H$ is injective on $S(U)$. Hence $y'' = 0$, and $x$~divides
$y$. By unique factorization in~$S(V^*)$ it follows that $L_n(V)$
divides $y$. So $\Ess(V)$ is the principal ideal generated
by~$L_n(V)$, and the free module on this one generator. Finally, the
definition of the essential ideal means that it is closed under the
action of the Steenrod algebra.
\end{proof}

\noindent
We finish off this section by recalling the action of the Steenrod
algebra on the cohomology of an elementary abelian $p$-group in the case
of an odd prime.
So let $p$ be an odd prime and $V$ an elementary abelian $p$-group.
Recall that the mod-$p$-cohomology ring is the free graded commutative algebra
\[
H^*(V,\f) \cong \f\lbrack x_1,\ldots,x_n\rbrack \otimes_{\f}
\Lambda(a_1,\ldots,a_n) \, ,
\]
where $a_i \in H^1(V,\f)$, $x_i \in H^2(V,\f)$, and $n$~is the rank
of~$V$. That is, $a_1,\ldots,a_n$ is a basis of the exterior copy
of~$V^*$, and $x_1,\ldots,x_n$ is a basis of the polynomial copy.
The product $a_1 a_2 \cdots a_n \in H^n(V,\f)$ is a basis of the top
exterior power $\Lambda^n(V^*)$. The Steenrod algebra~$\Ste$ acts on
the cohomology ring, making it an unstable $\Ste$-algebra with
$\beta(a_i)=x_i$ and $\rP^1(x_i) = x_i^p$. Observe that
$L_n(x_1,\ldots,x_n)$ is essential, for the same reason as in the
case $p=2$.

\section{The \Mui. invariants}
\label{section:Mui}
Let $k$ be a finite field and $V$ a finite dimensional $k$-vector space.
Consider the natural action of $\GL(V)$ on~$V^*$. The Dickson invariants
generate the invariants for the induced action of $\GL(V)$ on
the polynomial algebra~$S(V^*)$. But there is also an induced action
on the polynomial tensor exterior algebra $S(V^*) \otimes_k \Lambda(V^*)$,
and the \Mui. invariants are $\SL(V)$-invariants of this action:
see \Mui.'s original paper~\cite{Mui:ModInvts} as well as Crabb's modern
treatment~\cite{Crabb:DicksonMui}\@.

We shall need several properties of the \Mui. invariants.
For the convenience of the reader, we rederive these from scratch:
but see \Mui.'s papers \cite{Mui:ModInvts,Mui:CohoOps} and
Sum's work~\cite{Sum:Steenrod}\@.

\begin{notn}
Often we shall work with the direct sum decomposition
\[
H^*(V,\f) = \bigoplus_{r=0}^n N_r(V) \, ,
\]
where $n$~is the rank of~$V$ and we set
\[
N_r(V) = S(V^*) \otimes_{\f} \Lambda^r(V^*) \, .
\]
Observe that restriction to each subgroup respects this decomposition. This
means that the essential ideal is well-behaved with respect to this
decomposition:
\begin{equation}
\Ess(V) = \bigoplus_{r=0}^n N_r(V) \cap \Ess(V) \, .
\label{eqn:EssVdecomp}
\end{equation}
\end{notn}

\begin{defn}
Recall that $L_n(x_1,\ldots,x_n)$ is the determinant of the
$n \times n$-matrix
\[
C = \begin{pmatrix}
x_1 & x_2 & \cdots & x_n \\
\vdots & \vdots & \ddots & \vdots \\
x_1^{p^{n-1}} & x_2^{p^{n-1}} & \cdots & x_n^{p^{n-1}}
\end{pmatrix} \, ,
\]
where $C_{s,i} = x_i^{p^{s-1}}$ for $1 \leq s \leq n$.
For each such~$s$,
define $E(s)$ to be the matrix obtained from~$C$
by deleting row $s$ and then prefixing
$\begin{pmatrix} a_1 & a_2 & \cdots & a_n \end{pmatrix}$ as new first row:
so
\[
\det E(s) = \sum_{i=1}^n (-1)^{i+1} \gamma_{s,i} a_i \, ,
\]
where $\gamma_{s,i}$ is the
determinant of the minor of~$C$ obtained by removing row~$s$ and column~$i$.

Now define the \Mui. invariant $M_{n,s} \in H^*(V,\f)$
by $M_{n,s} = \det E(s)$. Note that our indexing differs from \Mui.'s:
our $M_{n,s}$ is his $M_{n,s-1}$.
\end{defn}

\begin{bsp}
So $M_{4,3} = \begin{vmatrix}
a_1 & a_2 & a_3 & a_4 \\
x_1 & x_2 & x_3 & x_4 \\
x_1^p & x_2^p & x_3^p & x_4^p \\
x_1^{p^3} & x_2^{p^3} & x_3^{p^3} & x_4^{p^3}
\end{vmatrix}$
and $\gamma_{2,3} = \begin{vmatrix}
x_1 & x_2 & x_4 \\ x_1^{p^2} & x_2^{p^2} & x_4^{p^2} \\
x_1^{p^3} & x_2^{p^3} & x_4^{p^3} \end{vmatrix}$.
\end{bsp}

\begin{lemma}
\label{lemma:Mns}
$M_{n,s} \in N_1(V) \cap \Ess(V)$.
\end{lemma}

\begin{proof}
By construction $M_{n,s} \in N_1(V)$.
Restricting to a maximal subgroup of~$V$ involves killing a nonzero
linear form on~$V^*$: that is, one imposes a linear dependence on
the~$a_i$ and consequently the same linear dependency on the~$x_i$.
So one obtains a linear dependency between the columns of~$E(s)$,
meaning that restriction kills $M_{n,s} = \det E(s)$.
\end{proof}

\begin{lemma}
\label{lemma:EssSquared}
$\Ess(V)^2 = L_n(V) \cdot \Ess(V)$.
\end{lemma}

\begin{proof}
As $L_n(V)$ is essential, the left hand side contains the right.
Now let $H$ be a maximal subgroup of~$V$. Then $H = \ker(a)$ for
some nonzero $a \in H^1(V,\f)$. Let $x = \beta(a) \in H^2$.
Observe that the kernel of restriction to~$H$ is generated by $a,x$.
Suppose that $f,g$ both lie in this kernel:
then we may write $f = f' a + f'' x$, $g = g' a + g'' x$, and so
$fg = (f'' g' \pm f' g'') a x + f'' g'' x^2$,
that is
$fg = x h$ for $h = (f''g' \pm f'g'')a + f''g'' x \in \ker \Res_H$.

Since $H^*(V,\f)$ is a free module over the unique factorization
ring $S(V^*)$, this means that $fg = L_n(V) \cdot y$ for some $y \in
H^*(V,\f)$. So $h = \frac{L_n(V)}{x} \cdot y$. As $\Res_H(h)=0$ and
$\Res_H\left(\frac{L_n(V)}{x}\right)$ is a non-zero divisor, we
deduce that $\Res_H(y)=0$. So $y \in \Ess(V)$.
\end{proof}

\begin{definition}
\label{defn:Mui}
Let $S = \{s_1,\ldots,s_r\} \subseteq \{1,\ldots,n\}$ be
a subset with $s_1 < s_2 < \cdots < s_r$. In view of
Lemmas \ref{lemma:Mns}~and \ref{lemma:EssSquared} we may define
the \Mui. invariant $M_{n,S} \in N_r(V) \cap \Ess(V)$
by
\[
M_{n,S} = \frac1{L_n(V)^{r-1}} M_{n,s_1} M_{n,s_2} \cdots M_{n,s_r} \, .
\]
Note in particular that $M_{n,\emptyset} = L_n(V)$.
\end{definition}

\begin{rk}
Observe that
\begin{equation}
M_{n,S} M_{n,T} = \begin{cases} \pm L_n(V) M_{n,S \cup T} &
\text{if $S \cap T = \emptyset$;} \\
0 & \text{otherwise.}
\end{cases}
\label{eqn:MnST}
\end{equation}
\end{rk}

\section{Joint annihilators}
\label{section:jointAnn}
\noindent
In this section we study the joint annihilators of the $M_{n,S}$ with
$\abs{S}=r$ as a means to prove Theorem~\ref{thm:free}\@.

\begin{lemma}
\label{lemma:jointAnn} The joint annihilator of
$M_{n,1},\ldots,M_{n,n}$ is $N_n(V)$.
\end{lemma}

\begin{proof}
The element $a_1 \ldots a_n$ is a basis for $\Lambda^n(V)$
and is clearly annihilated by each~$M_{n,s}$. Conversely, suppose that
$y\neq 0$ is annihilated by every~$M_{n,s}$. As
$M_{n,s} N_r(V) \subseteq N_{r+1}(V)$ we may assume without loss of generality
that $y \in N_r(V)$ for some~$r$.
Multiplying once or more by suitably chosen elements~$a_i$, we reduce
to the case $y \in N_{n-1}(V)$.

Denote by~$K$ the field of fractions of~$S(V^*)$, and let
$W = K \otimes_k \Lambda^{n-1}(V^*)$. Each $M_{n,s}$ induces
a linear form $\phi_s \colon W \rightarrow K$ given by
$\phi_s(w) a_1 \cdots a_n = M_{n,s} w$. By assumption, $y \neq 0$
lies in the kernel of every~$\phi_s$. A basis for~$W$ consists of
the elements $a_1 \cdots \widehat{a_r} \cdots a_n$ for $1 \leq r \leq n$,
where the hat denotes omission.
Now,
\[
M_{n,s} \cdot a_1 \cdots \widehat{a_r} \cdots a_n = (-1)^{r+1} \gamma_{s,r}
a_r \cdot a_1 \cdots \widehat{a_r} \cdots a_n \, ,
\]
and so
\[
\phi_s(a_1 \cdots \widehat{a_r} \cdots a_n) = \gamma_{s,r} \, .
\]
Now consider the matrix $\Gamma \in M_n(K)$ given by
$\Gamma_{s,r} = \gamma_{s,r}$. If one transposes and then
multiplies the $i$th row by $(-1)^i$ and the $j$th column by $(-1)^j$,
then one obtains the adjugate matrix of~$C$. As the determinant of~$C$ is
$L_n(V)$ and in particular nonzero, it follows that $\det \Gamma \neq 0$.

So by construction of~$\Gamma$, the $\phi_s$ form a basis of $W^*$.
So their common kernel is zero, contradicting our assumption on~$y$.
\end{proof}

\begin{corollary}
\label{coroll:jointAnn2}
The joint annihilator of $\{M_{n,S} \colon \abs{S} = r\}$ is
$\bigoplus_{s \geq n-r+1} N_s(V)$.
\end{corollary}

\begin{proof}
By induction on~$r$, Lemma~\ref{lemma:jointAnn} being the case $r=1$.
As $M_{n,S} \in N_{\abs{S}}(V)$ and $N_r(V) N_s(V) \subseteq N_{r+s}(V)$,
the annihilator is at least as large as claimed. Now suppose that
$y \in H^*(V,\f)$ does not lie in $\bigoplus_{s \geq n-r+1} N_s(V)$.
We may therefore write
\[
y = \sum_{s=0}^n y_s
\]
with $y_s \in N_s(V)$, and we know that $s_0 \leq n-r$ for $s_0 =
\min \{s \mid y_s \neq 0 \}$. As $y_{s_0} \neq 0$ and $y_{s_0} \not
\in N_n(V)$, Lemma~\ref{lemma:jointAnn} tells us that $y_{s_0}
M_{n,t} \neq 0$ for some $1 \leq t \leq n$.  As $y_{s_0} M_{n,t} \in
N_{s_0 + 1}(V)$, we conclude that $y M_{n,t}$ lies outside
$\bigoplus_{s \geq n-r+2} N_s(V)$. So the inductive hypothesis means
that there is some $T$~with $\abs{T} = r-1$ and $y M_{n,t} M_{n,T}
\neq 0$. So $y M_{n,S} \neq 0$ for $S = T \cup \{t\}$ and $\abs{S} =
r$: note that $t \in T$ is impossible.
\end{proof}

\begin{corollary}
\label{coroll:MnS}
Every $M_{n,S}$ is nonzero. For $S = \underline{n}
= \{1,\ldots,n\}$ we have
\[
\text{$M_{n,\underline{n}}$ is a nonzero scalar multiple of
$a_1 a_2 \cdots a_n$.}
\]
\end{corollary}

\begin{proof}
Observe that $M_{n,\underline{n}}$ is a scalar multiple of $a_1 \cdots a_n$
for degree reasons.
The case $r=n$ of Corollary~\ref{coroll:jointAnn2} says that $1 \in N_0(V)$
does not
annihilate $M_{n,\underline{n}}$ and therefore $M_{n,\underline{n}} \neq 0$.
But from Eqn~\eqref{eqn:MnST} we see that
every $M_{n,S}$ divides $L_n(V) M_{n,\underline{n}} \neq 0$.
\end{proof}

\begin{proof}[Proof of Theorem~\ref{thm:free}]
In view of Eqn~\eqref{eqn:EssVdecomp} it suffices to show that
for each~$r$ the \Mui. invariants $M_{n,S}$ with $\abs{S}=r$
are a basis of the $S(V^*)$-module $N_r(V) \cap \Ess(V)$. We observed in
Definition~\ref{defn:Mui} that these $M_{n,S}$ lie in this module.

So suppose that $y \in N_r(V) \cap \Ess(V)$. We should like there to be
$f_S \in S(V^*)$ such that
\begin{equation}
y = \sum_{\abs{S} = r} f_S M_{n,S} \, .
\label{eqn:pfFree}
\end{equation}
Note that for $T = \underline{n} - S$ we have
$M_{n,S} M_{n,T} = \pm L_n(V) M_{n,\underline{n}}$ by Eqn~\eqref{eqn:MnST}\@.
Define
$\eps_S \in \{+1,-1\}$ by $M_{n,S} M_{n,T} = \eps_S L_n(V) M_{n,\underline{n}}$.
So Eqn.~\eqref{eqn:pfFree} implies that we should define $f_S$ by
\[
f_S M_{n,\underline{n}} = \frac1{L_n(V)} \eps_S y M_{n,T} \, ,
\]
since $T \cap S' \neq \emptyset$ and therefore $M_{n,S'} M_{n,T}=0$
for all $S' \neq S$ with $\abs{S} = r$.
Note that this definition of~$f_S$ makes sense, as $y M_{n,T}$ lies
in both $N_r(V) N_{n-r}(V)=N_n(V)$ and $L_n(V) \Ess(V)$, the latter
inclusion coming from Lemma~\ref{lemma:EssSquared}\@.

With this definition of~$f_S$ we have
\[
\left(y - \sum_{\abs{S}=r} f_S M_{n,S}\right) M_{n,T} = 0
\]
for every $\abs{T} = n-r$. As $y -  \sum_{\abs{S}=r} f_S M_{n,S}$
lies in $N_r(V)$, this means that $y = \sum_{\abs{S}=r} f_S M_{n,S}$
by Corollary~\ref{coroll:jointAnn2}\@.

Finally we show linear independence. Suppose that $g_S \in S(V^*)$
are such that $\sum_{\abs{S} = r} g_S M_{n,S} = 0$. Pick one~$S$ and
set $T = \underline{n} - S$. Multiplying by $M_{n,T}$, we deduce
that $g_S=0$.
\end{proof}

\section{The action of the Steenrod algebra}
\label{section:Steenrod}
\noindent
To prepare for the proof of Theorem~\ref{thm:Steenrod} we shall study the
operation of the Steenrod algebra on the \Mui. invariants.

\begin{lemma}
\begin{xalignat}{2}
\beta (M_{n,s}) & = \begin{cases} L_n(V) & s = 1 \\ 0 & \text{otherwise}
\end{cases}  & \beta(L_n(V)) & = 0 \, .
\label{eqn:betaMns}
\end{xalignat}
For $0 \leq s \leq n-2$ we have:
\begin{xalignat}{2}
\rP^{p^s} (M_{n,r}) & = \begin{cases} M_{n,r-1} & r = s+2 \\ 0
& \text{otherwise}
\end{cases} & \rP^{p^s} (L_n(V)) & = 0 \, .
\label{eqn:rPMnS}
\end{xalignat}
\end{lemma}

\begin{proof}
One sees Eqn~\eqref{eqn:betaMns} by inspecting the determinants in the
definition of $M_{n,s}$ and $L_n(V)$.
The proof of  Eqn~\eqref{eqn:rPMnS} is also based on an inspection of these
determinants. Recall that $\rP^m(a_i)=0$ for every $m > 0$, and that
$\rP^m (x_i^{p^s})$ is zero too except for
$\rP^{p^s}(x_i^{p^s}) = x_i^{p^{s+1}}$.
We may use the Cartan formula
\[
\rP^m (xy) = \sum_{a+b=m} \rP^a(x)\rP^b(y)
\]
to distribute $\rP^{p^s}$ over the rows of the determinant. As $p^s$ cannot
be expressed as a sum of distinct smaller powers of~$p$, we only have
to consider summands where all of $\rP^{p^s}$ is applied to one row and
the other rows are unchanged. This will result in two rows being equal unless
it is the row consisting of the $x_i^{p^{s+1}}$ that is missing.
\end{proof}

\begin{lemma}
\label{lemma:SteenrodMnS}
Let $S = \{s_1,\ldots,s_r\}$ with $1 \leq s_1 < s_2 < \cdots < s_r \leq n$.
\begin{enumerate}
\item
\label{enum:betaMnS}
Suppose that $1 \not \in S$.
Then $M_{n,S} = \beta(M_{n,S\cup\{1\}})$.
\item
$L_n(V)^{r-1} \rP^m (M_{n,S}) = \rP^m(M_{n,s_1} \cdots M_{n,s_r})$
for each $m < p^{n-1}$.
\item
\label{enum:part3}
For $2 \leq u \leq n$ set $X = \{s \in S \mid s \leq u\}$ and
$Y = \{s \in S \mid s > u\}$. Then
\[
L_n(V) \rP^{p^{u-2}} (M_{n,S}) = \rP^{p^{u-2}}(M_{n,X}) \cdot M_{n,Y} \, .
\]
\item
For $1 \leq r \leq n$ and $0 < m < p^{n-1}$ one has
$\rP^m(M_{n,\{1,\ldots,r\}})=0$.
\item
\label{enum:part5}
For $2 \leq u \leq n$ one has
$\rP^{p^{u-2}} (M_{n,\{1,\ldots,u-2,u\}}) = M_{n,\{1,\ldots,u-1\}}$.
\end{enumerate}
\end{lemma}

\begin{proof}
Recall that
\begin{equation}
L_n(V)^r M_{n,S} = L_n(V) M_{n,s_1} \cdots M_{n,s_r} \, .
\label{eqn:recallMnS}
\end{equation}
The first two parts follow by applying
Equations \eqref{eqn:betaMns}~and \eqref{eqn:rPMnS}\@.

Recall that by the Adem relations each $\rP^m$ may be expressed in terms of
the $\rP^{p^s}$ with $p^s \leq m$.
So the third part follows from the second, since we deduce from
Eqn.~\eqref{eqn:rPMnS}
that $\rP^m (M_{n,s}) = 0$ if $0 < m \leq p^{u-2}$ and $s > u$.

Fourth part: By induction on~$r$. Follows for $r=1$ from
the Adem relations and Eqn~\eqref{eqn:rPMnS}\@. Inductive step:
Enough to consider $\rP^{p^s}$ for $0 \leq s \leq n-2$.
By the inductive hypothesis and a similar argument to the third
part, deduce that
\[
L_n(V) \rP^{p^s} (M_{n,\{1,\ldots,r\}}) = M_{n,\{1,\ldots,r-1\}}
\rP^{p^s} (M_{n,r}) \, .
\]
But this is zero by Eqn~\eqref{eqn:rPMnS}, since
$M_{n,\{1,\ldots,r-1\}} M_{n,r-1} = 0$.

Fifth part: Using the fourth part and an argument similar to the third
part, deduce that
\[
L_n(V) \rP^{p^{u-2}} (M_{n,\{1,\ldots,u-2,u\}})
= M_{n,\{1,\ldots,u-2\}} \rP^{p^{u-2}} (M_{n,u})
= M_{n,\{1,\ldots,u-2\}} M_{n,u-1} \, :
\]
but this is $L_n(V) M_{n,\{1,\ldots,u-1\}}$.
\end{proof}

\begin{proof}[Proof of Theorem~\ref{thm:Steenrod}]
We shall show that for every $M_{n,S}$ there is an element $\theta$ of
the Steenrod algebra with $M_{n,S} = \theta (M_{n,\underline{n}})$.
We do this by decreasing induction on $r = \abs{S}$. It is trivially
true for $r=n$, so assume now that $r < n$. Amongst the $S$~with
$\abs{S}=r$ we shall proceed by induction over~$u$, the smallest
element of $\underline{n} - S$. So
\[
S = \{1,\ldots,u-1\} \cup Y \quad \text{with $s > u$ for every $s \in Y$}\, .
\]
Part~\ref{enum:betaMnS} of Lemma~\ref{lemma:SteenrodMnS} covers the case
$u = 1$, so assume that $u \geq 2$.
Set $T = \{1,\ldots,u-2,u\}$. We complete the induction by showing that
$M_{n,S} = P^{p^{u-2}} (M_{n,T \cup Y})$.
Part~\ref{enum:part3} of Lemma~\ref{lemma:SteenrodMnS} gives us
\[
L_n(V) P^{p^{u-2}} (M_{n,T \cup Y}) = P^{p^{u-2}} (M_{n,T}) M_{n,Y} \, .
\]
But $P^{p^{u-2}} (M_{n,T}) = M_{n,\{1,\ldots,u-1\}}$,
by Part~\ref{enum:part5} of that lemma.
So $P^{p^{u-2}} (M_{n,T \cup Y}) = M_{n,S}$,
as claimed.
\end{proof}

\begin{rk}
Theorem~\ref{thm:free} shows that
the $S(V^*)$-module generated by the \Mui. invariants
$M_{n,S}$ is the essential ideal and therefore closed under the action of the
Steenrod algebra.  One may however see more directly that this $S(V^*)$-module
is Steenrod closed. This is observed for example in~\cite{Sum:Steenrod}\@.
In view of Lemma~ \ref{lemma:SteenrodMnS} and Equations \eqref{eqn:betaMns}~and
\eqref{eqn:rPMnS} it only remains to show that
$\rP^{p^{n-1}} (M_{n,s})$ lies in our $S(V^*)$-module.
Now $P^{p^{n-1}} (M_{n,n}) = 0$ by the unstable condition, so suppose
$s < n$. Recall that $M_{n,s}$ is a determinant, the last row
of the matrix having entries $x_i^{p^{n-1}}$. So applying
$P^{p^{n-1}}$ replaces these entries by $x_i^{p^n}$. But it is well known
that $x_i^{p^n}$ is an $S(V^*)$-linear combination of the
$x_i^{p^r}$ for $r \leq n-1$, and that the coefficients are independent of~$i$:
this is the ``fundamental equation'' in the sense of~\cite{Wilkerson:Dickson},
and the coefficients are the Dickson invariants $c_{n,r}$
in $S(V^*)$. Applying
$S(V^*)$-linearity of the determinant in the bottom row of the matrix,
one deduces that $P^{p^{n-1}} (M_{n,s})$ is an $S(V^*)$-linear combination
of the $M_{n,r}$.
\end{rk}

\bibliographystyle{abbrv}
\bibliography{../united}

\end{document}